\documentclass[12pt]{article}          
\usepackage{amssymb,amsgen,amsmath,amsbsy,amsopn,amsthm,amssymb}
\usepackage{latexsym}
\usepackage{comment}
\usepackage{graphicx}
\usepackage{curves}
\usepackage{psfrag}  
\usepackage{color}

\renewcommand{\refname}{References}

\newcommand{\K}{\mathcal{K}} 
\newcommand{\E}{\mathbb{E}^3}

\newtheorem{Th}{Theorem}

\title{Regular Polygonal Complexes of Higher Ranks in $\mathbb{E}^3$\\
\footnote{Appeared in {\em Model. Anal. Inform. Sist.\/} 20 (2013), 103--110. }}

\author{Egon Schulte\\
Northeastern University\\
Department of Mathematics\\
Boston, MA 02115, USA\\
schulte@neu.edu}

\date{\sl\small\today }

\begin{document}

\maketitle

\begin{abstract}
\noindent
The paper establishes that the rank of a regular polygonal complex in $\E$ cannot exceed~$4$, and that the only regular polygonal complexes of rank $4$ in $\E$ are the eight regular $4$-apeirotopes in $\E$. \\[.1in]
\end{abstract}

\section{Introduction} 
A geometric polygonal complex $\K$ of rank $n$ in Euclidean 3-space $\E$ is a discrete incidence structure made up from finite or infinite, planar or skew polygons, assembled in a careful fashion into families of geometric polygonal complexes of smaller ranks. As combinatorial objects they are incidence complexes of rank $n$ with polygons as $2$-faces, that is, abstract polygonal complexes of rank $n$ (see \cite{ds,es1}). A geometric polygonal complex $\K$ is {\em regular\/} if $\K$ has a flag-transitive geometric symmetry group.

The regular polygonal complexes of rank $3$ in $\E$ are completely known. They comprise fourty-eight regular polyhedra as well as twenty-five regular complexes which are not polyhedra. The regular polyhedra were discovered by Gr\"unbaum~\cite{gr1} and Dress~\cite{d1,d2}, and are described in detail in McMullen \& Schulte~\cite[Ch. 7E]{arp} and McMullen~\cite{grp}. The regular complexes which are not polyhedra were recently classified in Pellicer \& Schulte~\cite{pelsch1,pelsch2}.

The present paper proves that a regular polygonal complex in $\E$ cannot have rank larger than $4$, and that the only regular polygonal complexes of rank $4$ in $\E$ are the eight regular $4$-apeirotopes in $\E$ described in McMullen \& Schulte~\cite[Ch. 7F]{arp}. 

\section{Incidence complexes}
\label{inc}

Following \cite{ds,es1}, an {\em incidence complex of rank\/} $n$, or simply an {\em $n$-complex}, is a partially ordered set $\mathcal{K}$ with a strictly monotone rank function with range $\{-1,0,\ldots,n\}$ satisfying the following conditions. The elements of rank $j$ are called the {\em $j$-faces\/} of $\mathcal{K}$, or {\em vertices}, {\em edges\/} and {\em facets\/} of $\mathcal{K}$ if $j=0$, $1$ or $n-1$, respectively. Two faces $F$ and $G$ are said to be {\em incident\/} if $F\leq G$ or $G\leq F$. Each {\em flag\/} (maximal totally ordered subset) of $\mathcal{K}$ is required to contain exactly $n+2$ faces, including a unique minimal face $F_{-1}$ (of rank $-1$) and a unique maximal face $F_{n}$ (of rank $n$) as {\em improper\/} faces. We often find it convenient to suppress the improper faces in designating flags. Two flags $\Phi$ and $\Psi$ of $\K$ are called {\em adjacent\/} if they differ in exactly one face; if this face is an $j$-face for some~$j$ (with $0\leq j\leq n-1$), then $\Phi$ and $\Psi$ of $\K$ are {\em $j$-adjacent\/}. Further, we ask that $\mathcal{K}$ be {\em strongly flag-connected\/}, meaning that any two flags $\Phi$ and $\Psi$ of $\mathcal{P}$ can be joined by a sequence of flags $\Phi=\Phi_{0},\Phi_{1},\ldots,\Phi_{l-1},\Phi_{l}=\Psi$, all containing $\Phi\cap\Psi$, such that $\Phi_{i-1}$ and $\Phi_{i}$ are {\em adjacent\/} (differ by exactly one face) for each $i$. Finally, $\K$ has the property that if $F$ and $G$ are incident faces of ranks $j-1$ and $j+1$ for some $j$, then there are {\em at least two\/} $j$-faces $H$ such that $F<H<G$. 

With regards to this latter condition, note that the corresponding condition required in \cite{ds,es1} is more restrictive, in that, for each $j$, a fixed number $k_j$ of $j$-faces $H$ is required to lie between $F$ and $G$. This stronger condition is always satisfied for highly symmetric complexes like those studied in this paper. 

When $F$ and $G$ are two faces of a complex $\mathcal{K}$ with $F \leq G$, we call 
\[ G/F := \{H \mid F \leq H \leq G\}\] 
a {\em section\/} of $\mathcal{K}$. We usually identify a face $F$ with the section $F/F_{-1}$. For a face $F$ we also call $F_{n}/F$ the {\em co-face of $\mathcal{K}$ at\/} $F$, or the {\em vertex-figure at\/} $F$ if $F$ is a vertex.  

An incidence complex $\K$ is said to be {\em regular\/} if its ({\em combinatorial\/}) {\em automorphism group\/} $\Gamma(\mathcal{K})$ is transitive on the flags. 

A complex $\K$ is an ({\em abstract}) {\em polytope\/} if, for all $j$ and all $(j-1)$-faces $F$ and $(j+1)$-faces $G$ with $F<G$, there are exactly {\em two\/} $j$-faces $H$ such that $F<H<G$. This class of complexes has attracted a lot of attention (see \cite{arp}). 

More generally, an incidence complex $\K$ is called an {\em abstract polygonal complex\/} if each $2$-face of $\K$ is isomorphic to (the face lattice of) a convex polygon or an (infinite) apeirogon.

For $0\leq k\leq n-1$, the {\em $k$-skeleton\/} of a complex $\K$ of rank $n$ is the incidence complex of rank $k+1$, whose faces comprise the $n$-face of $\K$ and all faces of $\K$ of rank less than or equal to $k$, with the partial order inherited from $\K$.  

\section{Geometric polygonal complexes}
\label{real}

Our definition of a ``realization'' of an abstract polygonal complex is inspired by the corresponding notion for abstract polytopes (see \cite{mr,arp}). 

A (faithful) {\em realization\/} of an abstract polygonal complex $\K$, again denoted by $\K$, is derived inductively from a bijection $\beta$ of the vertex-set of $\K$ into some Euclidean space $\mathbb{E}$. The vertices of $\K$ are mapped by $\beta$ to the {\em vertices\/} of the realization. Then each $1$-face of $\K$ can be viewed as being mapped under $\beta$ to a line segment, called an {\em edge} of the realization, joining the images of the vertices of the $1$-face under $\beta$. Moving up in rank, each $2$-face of $\K$, which by assumption is isomorphic to a convex polygon or an apeirogon, is mapped to the finite or infinite polygon in $\mathbb{E}$, a {\em $2$-face\/} of the realization, made up from the edges of the realization that are the images under $\beta$ of the $1$-faces of the given $2$-face of $\K$. More generally, from the $j$-faces of the realization we then obtain each {\em $(j+1)$-face\/} of the realization as a family of $j$-faces of the realization, namely those corresponding under $\beta$ to the $j$-faces of the underlying $(j+1)$-face of $\K$. Finally, then, when $j=n-1$ we arrive at the desired realization of $\K$ determined by $\beta$. We also call the resulting structure a {\em geometric polygonal complex\/}.

Clearly, in order for this process to produce a faithful geometric copy of the given abstract polygonal complex $\K$ we must assume that, for each $j$, each $(j+1)$-face of $\K$ is uniquely determined by its $j$-faces. Throughout we will be working under this assumption. Note, however, that an abstract complex may not satisfy this assumption; the regular map $\{3,6\}_{(2,0)}$ on the $2$-torus is an example.

Alternatively we can proceed more directly and define geometric polygonal complexes with less explicit reference to realizations of abstract polygonal complexes, at least for small ranks, beginning with rank $2$. Here we restrict ourselves to Euclidean $3$-space $\E$, although similar notions carry over to realizations in higher-dimensional Euclidean spaces.

Following \cite{gr1}, a {\em finite polygon\/}, or an {\em $n$-gon\/}, in $\mathbb{E}^3$ consists of a sequence $(v_1, v_2, \dots, v_n)$ of $n$ distinct points, as well as of the line segments $(v_1, v_2), (v_2,v_3), \ldots, (v_{n-1},v_n),(v_n, v_1)$.  A (discrete) {\em infinite polygon\/}, or {\em apeirogon\/}, in $\mathbb{E}^3$ similarly consists of an infinite sequence of distinct points $(\dots, v_{-2},v_{-1}, v_0, v_1, v_2, \dots)$, as well as of the line segments $(v_i, v_{i+1})$ for each $i$, such that each compact subset of $\mathbb{E}^3$ meets only finitely many line segments.  In either case the points are the {\em vertices\/} and the line segments the {\em edges\/} of the polygon.

Then following \cite{pelsch1}, a {\em geometric polygonal complex $\K$ of rank $3$}, or simply a {\em geometric polygonal 3-complex}, in $\mathbb{E}^{3}$ is a triple $(\mathcal{V},\mathcal{E},\mathcal{F}$) consisting of a set $\mathcal{V}$ of points, called {\em vertices}, a set $\mathcal{E}$ of line segments, called {\em edges}, and a set $\mathcal{F}$ of (finite or infinite) polygons, called {\em faces} (here to mean {\em $2$-faces\/}), satisfying the following properties.\\[-.25in]
\begin{itemize}
\item[(a)] The graph $(\mathcal{V},\mathcal{E})$, the {\em edge graph\/} or {\em net\/} of $\K$, is connected.\\[-.3in]
\item[(b)] The vertex-figure of $\K$ at each vertex of $\K$ is connected. By the {\em vertex-figure\/} of $\K$ at a vertex $v$ we mean the graph, possibly with multiple edges, whose vertices are the vertices of $\K$ adjacent to $v$ and whose edges are the line segments $(u,w)$, where $(u, v)$ and $(v, w)$ are edges of a common face of $\K$. (There may be more than one such face in $\K$, in which case the corresponding edge $(u,w)$ of the vertex-figure at $v$ has multiplicity given by the number of such faces.)\\[-.3in]
\item[(c)] Each edge of $\K$ is contained in at least two faces of $\K$. (This requirement is less restrictive than the corresponding condition given in~\cite{pelsch1}. For highly symmetric complexes like those discussed here, the two conditions are equivalent.) \\[-.3in]
\item[(d)] $\K$ is {\em discrete\/}, in the sense that each compact subset of $\mathbb{E}^{3}$ meets only finitely many faces of $\K$.
\end{itemize}

Each geometric polygonal complex of rank $3$ in $\E$ gives an incidence complex of the same rank. In fact, a quicker definition would consist of saying that $\K$ as above consists of a triple $(\mathcal{V},\mathcal{E},\mathcal{F}$), which, when ordered by inclusion, gives an abstract polygonal complex of rank $3$. 

Proceeding with rank $4$ structures, a {\em geometric polygonal complex $\K$ of rank $4$\/}, or simply a {\em geometric polygonal 4-complex}, in $\mathbb{E}^{3}$ is a $4$-tuple $(\mathcal{F}^0,\mathcal{F}^1,\mathcal{F}^2,\mathcal{F}^3)$ consisting of a set $\mathcal{V}=\mathcal{F}^0$ of points, called {\em vertices}, a set $\mathcal{E}=\mathcal{F}^1$ of line segments, called {\em edges}, a set $\mathcal{F}=\mathcal{F}^2$ of (finite or infinite) polygons, called {\em $2$-faces\/}, and a set $\mathcal{F}^3$ of geometric polygonal $3$-complexes, called {\em $3$-faces\/}, such that~$\K$, partially ordered by inclusion, gives an abstract polygonal  complex of rank $4$ (with the same rank function).

We can move up further in rank and similarly define higher rank geometric polygonal complexes in $\E$. The basic set-up is exactly the same. The $j$-faces of an $n$-complex are built from the $(j-1)$-complexes representing the $(j-1)$-faces of the $j$-face.

A geometric polygonal complex is a {\em geometric polytope\/}, or a {\em geometric polyhedron\/} when the rank is $3$, if the underlying incidence complex is an abstract polytope. An {\em apeirotope\/} or {\em apeirohedron\/} is an infinite polytope or polyhedron, respectively. 

A geometric polygonal complex $\K$ is called ({\em geometrically\/}) {\em regular\/} if its ({\em geometric\/}) {\em symmetry group\/} $G(\mathcal{K})$ in $\E$ is transitive on the flags of $\K$. The symmetry group $G(\K)$ of a polygonal complex $\K$ is a (generally proper) subgroup of its combinatorial automorphism group $\Gamma(\K)$. In fact, $G(\K)$ is a flag-transitive subgroup of $\Gamma(\K)$ if $\K$ is regular. For a geometrically regular polytope, the symmetry group coincides with the full automorphism group and is necessarily simply flag-transitive. However, this is not true in general for geometrically regular polygonal complexes.

We later require the classification of the regular apeirotopes of rank $4$ in $\E$ (see~\cite[Ch. 7F]{arp}). There are eight examples shown in the following display:
\[\begin{matrix}
\{4,3,4\} = \{\{4,3\},\{3,4\}\}, & \{\{4,6 \mid 4\},\{6,4\}_{3}\};\\[.02in]
\{\{\infty,3\}_{6} \# {\{\,\}},\{3,3\}\}, & \{\{\infty,4\}_{4} \#
\{\infty\},\{4,3\}_{3}\};  \\[.02in]
\{\{\infty,3\}_{6} \# {\{\,\}},\{3,4\}\}, & \{\{\infty,6\}_{3} \#
\{\infty\},\{6,4\}_{3}\};  \\[.02in]
\{\{\infty,4\}_{4} \# {\{\,\}},\{4,3\}\}, & \{\{\infty,6\}_{3} \#
\{\infty\},\{6,3\}_{4}\}. 
\end{matrix} \] 
For notation we refer to \cite[Ch. 7F]{arp}. The eight examples fall into four pairs of ``Petrie-duals" (listed in the same row), where the apeirotopes in each pair share the same $2$-skeleton. The familiar cubical tessellation $\{4,3,4\}$ and its Petrie-dual listed in the first row both have (finite) square faces. All other apeirotopes have (infinite planar) zigzag faces.  

\section{Groups of regular incidence complexes}
\label{group}

The structure of a regular incidence complex $\K$ can be completely described in terms of a well-behaved system of generating subgroups of any flag-transitive subgroup $\Lambda$ of the full automorphism group $\Gamma(\K)$ of $\K$ (see \cite{es1}). 

Let $\K$ be a regular incidence complex of positive rank $n$, and let $\Lambda$ be any flag-transitive subgroup of $\Gamma(\K)$. Let $\Phi:=\{F_0,\ldots,F_{n-1}\}$ be a (fixed) {\em base flag\/} of $\K$, where $F_i$ denotes the $i$-face in $\Phi$ for each $i$. For $\Omega\subseteq \Phi$, let $\Lambda_\Omega$ denote the stabilizer of $\Omega$ in $\Lambda$, so in particular $\Lambda_\Phi$ is the stabilizer of $\Phi$ itself, and $\Lambda_\emptyset = \Lambda$. For $i=-1,0,\ldots,n$ define the subgroup 
\[ R_{i} :=  \Lambda_{\Phi\setminus\{F_i\}} = \langle \varphi\in \Lambda \mid F_j\varphi =F_j \mbox{ for all } j\neq i\rangle .\]
Then $R_{-1}=R_{n}=\Lambda_\Phi \leq R_i$ for each $i=0,\ldots,n-1$, and $|R_{i}:R_{-1}|-1$ is the number of flags $i$-adjacent to $\Phi$. Moreover, 
\begin{equation}
\label{comm}
R_i \cdot R_j = R_j \cdot R_i \qquad (-1\leq i <j-1\leq n-1), \
\end{equation}
as products of subgroups. 

The subgroup $\Lambda$ is generated by the {\em distinguished generating subgroups\/} $R_0,\ldots,R_{n-1}$, that is, 
\[ \Lambda = \langle R_0,\ldots,R_{n-1} \rangle. \]
In fact, a stronger statement holds. When $\emptyset\neq I\subseteq\{-1,0,\ldots,n\}$ set $\Lambda_{I} := \langle R_i \mid i\in I\rangle$; and when $I=\emptyset$ set $\Lambda_I := R_{-1}=\Lambda_\Phi$. Note that $\Lambda_{I}=\Lambda_{I\setminus\{-1,n\}}$ for each subset $I$. Then, for each $\Omega\subseteq\Phi$, 
\[ \Lambda_\Omega = \langle R_i \mid F_{i}\not\in\Omega\rangle = \Lambda_{\{i\mid F_i \not\in\Omega\}}.\]
Here both sides coincide with $\Lambda_\Phi$ if $\Omega=\Phi$. In addition, the following {\em intersection property\/} holds: 
\begin{equation}
\label{intsec}
\Lambda_I \cap \Lambda_J = \Lambda_{I\cap J}\qquad (I,J\subseteq \{-1,0,\ldots,n\}) . 
\end{equation}
The partial order on $\K$ can be completely described in terms of the subgroups $R_{-1},R_0,\ldots,R_n$ of $\Lambda$. In fact, 
\[ F_{i}\varphi \leq F_{j}\psi\; \leftrightarrow\; \psi^{-1}\varphi \in \Lambda_{\{i+1,\ldots,n\}}\Lambda_{\{-1,0,\ldots,j-1\}}
\;\quad (-1\leq i\leq j\leq n;\, \varphi,\psi\in\Lambda) ,\]
or equivalently,
\[ F_{i}\varphi \leq F_{j}\psi\, \leftrightarrow \,
\Lambda_{\{-1,0,\ldots,n\}\setminus\{i\}}\varphi \cap \Lambda_{\{-1,0,\ldots,n\}\setminus\{j\}}\psi \neq \emptyset
\quad (-1\leq i\leq j\leq n;\, \varphi,\psi\in\Lambda). \]

Conversely, if a group $\Lambda$ has a system of generating subgroups $R_{-1},R_0,\ldots,R_n$ with properties (\ref{comm}) and (\ref{intsec}), then $\Lambda$ is a flag-transitive subgroup of the full automorphism group of a regular incidence complex $\K$ of rank $n$ (see \cite{es1}).

For abstract regular polytopes, the groups $R_0,\ldots,R_{n-1}$ are generated by involutions $\rho_0,\ldots,\rho_{n-1}$, where $\rho_i$ maps the base flag to its unique $i$-adjacent flag (see \cite[Ch. 2]{arp}). Also, $\Lambda_\Phi$ is trivial in this case.

\section{Rank $3$}

The regular polygonal complexes of rank $3$ in $\E$ have been completely enumerated. They comprise fourty-eight regular polyhedra and twenty-five regular complexes which are not polyhedra.

The regular polyhedra in $\E$ were discovered by Gr\"unbaum~\cite{gr1} and Dress~\cite{d1,d2} and are sometimes referred to as the Gr\"unbaum-Dress polyhedra. For a quicker method of arriving at the classification see also McMullen \& Schulte~\cite{ms} and \cite[Ch. 7E]{arp}. The regular polyhedra include many well-known figures such as the Platonic solids, the Kepler-Poinsot polyhedra and the Petrie-Coxeter polyhedra (see \cite{cox1,cox2}). 

The regular polygonal complexes in $\E$ which are not polyhedra were enumerated in Pellicer \& Schulte~\cite{pelsch1,pelsch2}. The flag-stabilizers of their symmetry groups are either trivial or have order $2$, and the methods of enumeration in these two cases are quite different. 

There are four regular polygonal complexes with a non-trivial flag-stabilizer, and each is the $2$-skeleton of a regular $4$-apeirotope in $\E$. The eight regular $4$-apeirotopes in $\E$ fall into four pairs of Petrie-duals, and the apeirotopes in each pair have the same $2$-skeleton (see~\cite[Ch. 7F]{arp}). Thus each regular polygonal complex with a non-trivial flag-stabilizer is the common $2$-skeleton of two regular $4$-apeirotopes. 

There are also twenty-one regular polygonal complexes with a trivial flag-stabilizers, which are not polyhedra; accordingly, these complexes are referred to as the simply flag-transitive regular polygonal complexes. 

\section{Higher ranks}

In this section we bound the rank of a regular polygonal complex in $\E$ by~$4$, and show that the eight regular $4$-apeirotopes mentioned earlier are the only examples in rank $4$. 

Suppose $\K$ is a regular polygonal complex of rank $n\geq 4$. Then its geometric symmetry group $G(\K)$ is a flag-transitive subgroup of $\Gamma(\K)$ to which the above structure results for flag-transitive subgroups of $\Gamma(\K)$ apply with $\Lambda=G(\K)$. In particular, $G(\K)=\langle R_0,\ldots,R_{n-1}\rangle$, where $R_0,\ldots,R_{n-1}$ are the distinguished generating subgroups of $G(\K)$ defined with respect to a base flag $\Phi=\{F_0,\ldots,F_{n-1}\}$ of $\K$. Since $\K$ is a polygonal complex, the base $2$-face $F_2$ is a regular polygon in $\E$, planar or non-planar. 

Now each of the subgroups $R_{3},\ldots,R_{n-1}$ consists of isometries stabilizing $F_0,F_1,F_2$ and hence acting trivially (pointwise) on the entire affine hull of the polygon $F_2$. This forces $F_2$ to be planar. In fact, otherwise the subgroups $R_{3},\ldots,R_{n-1}$ would have to be trivial; but this is impossible since $n\geq 4$. On the other hand, if $F_2$ is planar, then a nontrivial element from $R_{3},\ldots,R_{n-1}$ could only be the euclidean reflection, $\rho_3$ (say), in the affine hull of $F_2$. In this case we can immediately exclude the possibility that $n>4$; in fact, if $n>4$ then necessarily $R_{3}=R_{4}=\langle\rho_3\rangle$, which is impossible. This then only leaves the possibility that $n=4$, the face $F_2$ is planar, $R_{3}=\langle\rho_3\rangle$, and there are just two facets of $\K$ meeting at $F_2$. 

Now suppose that $n=4$ and that $F_2$ and $\rho_3$ are as described. Then the $2$-skeleton $\mathcal{L}$ of $\K$ is a regular polygonal complex of rank $3$ in $\E$ whose symmetry group contains $G(\K)$. In particular, $\mathcal{L}$ has {\em face mirrors\/}, meaning that $\mathcal{L}$ has planar $2$-faces and that the affine hulls of the $2$-faces are mirrors of plane symmetries of $\mathcal{L}$. In fact, $\rho_3$ is a reflective symmetry of $\mathcal{L}$ in the affine hull of the $2$-face $F_2$ of $\mathcal{L}$, and its conjugates under $G(\K)$ (or $G(\mathcal{L})$) provide all reflective symmetries in affine hulls of $2$-faces of $\mathcal{L}$. In particular, $\mathcal{L}$ is not simply flag-transitive, since $\rho_3$ stabilizes the flag $\{F_0,F_1,F_2\}$ of $\mathcal{L}$.

We now appeal to the classification in \cite{pelsch1,pelsch2} of the regular polygonal complexes of rank $3$ in $\E$. In fact, it was shown in \cite{pelsch1} that a regular polygonal $3$-complex with face mirrors in $\E$ is necessarily the $2$-skeleton of a regular $4$-apeirotope $\mathcal{P}$ in $\E$ with the same symmetry group, and that the fourth distinguished generator of the symmetry group $G(\mathcal{P})$ of $\mathcal{P}$ is given by the reflection in the planar base $2$-face. Thus $\mathcal{L}$ is the $2$-skeleton of a regular $4$-apeirotope $\mathcal{P}$ with $G(\mathcal{P})=G(\mathcal{L})$. It is also known that the number of $2$-faces, $r$, on an edge of $\mathcal{P}$ and thus of $\mathcal{L}$ must be $3$ or $4$.

Now the facets of $\K$ are also regular polygonal complexes of rank $3$ in $\E$ and have their $2$-faces among the $2$-faces of $\mathcal{L}$. In particular, the basic facet $F_3$ of $\K$ gives rise to the complex $\mathcal{K}_{3}:=F_{3}/F_{-1}$ whose full symmetry group $G(\mathcal{K}_3)$ contains the stabilizer $\langle R_{0},R_{1},R_{2}\rangle$ of $F_3$ in $G(\K)$ as a flag-transitive subgroup. We must show that $\mathcal{K}_3$ is a polyhedron. Then, since there are just two facets of $\K$ meeting at $F_2$, the complex $\K$ itself would have to be a $4$-apeirotope. 

Since $\K_3$ is a subcomplex of $\mathcal{L}$, the number of $2$-faces, $r_3$, on an edge of $\K_3$ must be $2$, $3$ or $4$. We wish to show that $r_{3}=2$. Clearly, we cannot have $r_{3}=r$ since $\mathcal{L}$ is connected. In fact, otherwise, every $2$-face of $\K$ containing an edge of $\K_3$ would also have to be a $2$-face of $\K_3$; but then connectedness would force all $2$-faces of $\K$ to be $2$-faces of $\K_3$, which is impossible. Hence it remains to exclude the case when $r=4$ and $r_{3}=3$. Now suppose $r=4$ and $r_{3}=3$. Then we know from \cite{pelsch1} that the pointwise stabilizer of the base edge $F_{1}$ of $\K_3$ in $G(\K_3)$ is a cyclic group $C_3$ or a dihedral group $D_3$. Either way, the affine hulls of the (planar) faces of $\K_3$ meeting at $F_1$ must be inclined at $120^\circ$. On the other hand, the affine hulls of the four $2$-faces of $\mathcal{P}$ (or $\K$ or $\mathcal{L}$) meeting at $F_1$ either coincide (when the faces are opposite relative to $F_1$) or are inclined at $90^\circ$. In any case, the two scenarios are incompatible and rule out the possibility that $r=4$ and $r_{3}=3$. Thus $r_{3}=2$. 

At this point we know that both $\K$ and $\mathcal{P}$ are regular $4$-apeirotopes with a common $2$-skeleton. In particular, $\K$ must be among the eight regular $4$-apeirotopes in $\E$. Another appeal to \cite{pelsch1} then shows that $\K$ must in fact coincide with either $\mathcal{P}$ itself or with the Petrie-dual of $\mathcal{P}$.

In summary we have established the following theorem.

\begin{Th}
\label{mainthm}
There are no regular geometric polygonal complex in $\E$ of rank $n\geq 5$. The only regular geometric polygonal complexes of rank $4$ in $\E$ are the eight regular $4$-apeirotopes in $\E$. 
\end{Th}

\renewcommand{\refname}{References}

\end{document}